\documentclass[12pt]{amsart}
\usepackage[dvips]{graphicx}
\usepackage{amsmath,graphics}
\usepackage{xypic}
\usepackage{amsfonts,amssymb}
\theoremstyle{plain}
\newtheorem{theorem*}{Theorem}
\newtheorem*{lemma*}{Lemma}
\newtheorem{corollary*}{Corollary}
\newtheorem*{proposition*}{Proposition}
\newtheorem{conjecture*}{Conjecture}
\newtheorem{theorem}{Theorem}[section]
\newtheorem{lemma}[theorem]{Lemma}
\newtheorem{corollary}[theorem]{Corollary}
\newtheorem{proposition}[theorem]{Proposition}

\theoremstyle{remark}

\newtheorem*{remark}{Remark}

\theoremstyle{definition}
\newtheorem{defn}[theorem]{Definition}

 \def\Q{\Bbb{Q}}  \def\Z{\Bbb{Z}} \def\R{\Bbb{R}} 
    
 \def\a{\alpha}   \def\bp{\begin{pmatrix}}
\def\sm{\setminus} \def\ep{\end{pmatrix}} \def\bn{\begin{enumerate}} 
 \def\rank{\mbox{rank}} \def\div{\mbox{div}} \def\en{\end{enumerate}}
\def\ba{\begin{array}} \def\ea{\end{array}}  
 \def\S{\Sigma}  \def\a{\alpha} \def\b{\beta} \def\ti{\tilde}
  \def\im{\mbox{Im}} 
  
\def\ker{\mbox{Ker}}\def\be{\begin{equation}} \def\ee{\end{equation}}

   \def\zgt{\Z[G][t^{\pm 1}]}

\def\zt{\Z[t^{\pm 1}]}    
\def\w{\omega}   
    \def\fr12{\frac{1}{2}} \def\z12{\Z[\fr12]}

\def\tpm {[t^{\pm 1}]}

\def\i{\iota}

\begin{document}
\title[Symplectic $S^{1} \times N^3$]{Symplectic $\bf S^{1} \times N^3$, subgroup separability, and vanishing Thurston norm}
\author{Stefan Friedl}
\address{Universit\'e du Qu\'ebec \`a Montr\'eal\\ Montr\'eal, Qc H3C 3P8\\ Canada}
\email{sfriedl@gmail.com}
\author{Stefano Vidussi}
\address{Department of Mathematics\\ University of California\\
Riverside, CA 92521\\ USA} \email{svidussi@math.ucr.edu} \thanks{S.
Vidussi was partially supported by NSF grant \#0629956.}
\date{May 1, 2007}
\begin{abstract} Let $N$  be a closed, oriented $3$--manifold.
A folklore conjecture states that $S^{1} \times N$ admits a
symplectic structure if and only if $N$ admits a fibration over the
circle. We will prove this conjecture in the case when $N$ is
irreducible and its fundamental group satisfies appropriate
subgroup separability conditions. This statement includes
$3$--manifolds with vanishing Thurston norm, graph
manifolds and $3$--manifolds with surface subgroup separability (a
condition satisfied conjecturally by all hyperbolic $3$--manifolds).
Our result covers, in particular, the case of $0$--framed surgeries
along knots of genus one. The statement follows from the proof that
twisted Alexander polynomials decide fiberability for all the
$3$--manifolds listed above. As a corollary, it follows that twisted
Alexander polynomials decide if a knot of genus one is fibered.
\end{abstract}
\dedicatory{Dedicated to the memory of Xiao-Song Lin}
\maketitle
%\vspace{0.5cm} \hfill{\it Dedicated to the memory of Xiao-Song Lin}
%\\[0.5cm]

\section{Introduction}

Let $N$ be a 3--manifold. Throughout the paper, unless otherwise
stated, we will assume that all 3--manifolds are closed, oriented
and connected, all surfaces are oriented, and all homology and
cohomology groups have integer coefficients.

Thurston \cite{Th76} showed that if $N$ admits a fibration over
$S^{1}$, then $S^{1} \times N$ is symplectic, i.e. it can be endowed
with a closed, non--degenerate $2$--form $\w$.

It is natural to ask whether the converse of this statement holds
true. Interest on this question was motivated by Taubes' results in
the study of symplectic $4$-manifolds (see \cite{Ta94,Ta95}), that
gave some initial evidence to an affirmative answer to this
question. We can state this problem in the following form:

\begin{conjecture*} \label{conjfolk}
Let $N$ be a 3--manifold. If $S^1\times N$ is symplectic, then there
exists $\phi\in H^1(N)$ such that $(N,\phi)$ fibers over $S^1$.
\end{conjecture*}

Here we say that $(N,\phi)$ fibers over $S^{1}$ if the homotopy
class of maps $N \to S^1$ determined by $\phi \in H^{1}(N) =
[N,S^{1}]$ contains a representative that is a fiber bundle over
$S^{1}$; in that case, we will also say that $\phi$ is a
\textit{fibered class}.

Assuming the Geometrization Conjecture it is shown in \cite{McC01} that if $S^1\times N$ is symplectic,
then $N$ is prime. Excluding the trivial case $N=S^1\times S^2$ we will therefore restrict ourselves to studying Conjecture
\ref{conjfolk} for irreducible $3$--manifolds.

In \cite{FV06}, we suggested an approach to Conjecture
\ref{conjfolk} based on the study of twisted Alexander polynomials
$\Delta_{N,\phi}^\a$ of $N$ associated to some $\phi \in H^1(N)$ and
an epimorphism $\a$ of the fundamental group of $N$ onto a finite
group $G$. This approach, while relying on results from
Seiberg-Witten theory and symplectic topology, embeds Conjecture
\ref{conjfolk} in questions related with group theory for
$3$--manifold groups and the theory of covering spaces.

Precisely, in \cite{FV06} we showed that Conjecture \ref{conjfolk} is implied by the following (perhaps stronger) conjecture:

\begin{conjecture*} \label{conjcha}
Let $N$ be a $3$--manifold and let $\phi \in H^{1}(N)$ be a
primitive class such that for \textit{any} epimorphism onto a finite
group $\a: \pi_{1}(N) \rightarrow G$ the twisted Alexander
polynomial $\Delta^{\a}_{N,\phi} \in \Z[t^{\pm 1}]$ is monic and
$\deg\Delta_{N,\phi}^{\a} = |G| \, \|\phi\|_{T} + 2 \div \,
\phi_{G}$. Then $(N,\phi)$ fibers over $S^{1}$. \end{conjecture*}

(Notation and definitions relevant to this conjecture are presented
in Sections \ref{sw} and \ref{sectionufd}.)

Specifically, Theorem 4.3 of \cite{FV06} asserts that the conditions on the twisted Alexander polynomial $\Delta_{N,\phi}^{\a}$ required in Conjecture \ref{conjcha} are satisfied by the K\"unneth component in $H^1(N)$ of the class $[\omega] \in H^2(S^1 \times N)$ of an integral symplectic form on $S^1 \times N$.

In this paper we will collect some dividends from this approach.
Precisely, we will show in Theorem \ref{mainthm} that Conjecture
\ref{conjcha} holds true for a class of $3$--manifolds whose
fundamental group satisfies appropriate subgroup separability
conditions. For sake of exposition, we will quote here a slightly
weaker version of Theorem \ref{mainthm}. In order to state it,
recall that a subgroup $A \subset \pi_1(N)$ of the fundamental group
of a $3$-manifold is \textit{separable} if for any $g \in \pi_1(N)
\setminus A$ there exists an epimorphism $\a : \pi_1(N) \to G$,
where $G$ is a finite group, such that $\a(g) \notin \a(A)$.
We also say that $\pi_1(N)$ is \textit{surface subgroup separable} if any surface group  $A\subset \pi_1(N)$  is separable.

We have the following:

\begin{theorem*} \label{whatever}
Let $N$ be an irreducible $3$--manifold and let $\phi \in H^{1}(N)$
be a primitive class such that for \textit{any} epimorphism onto a
finite group $\a: \pi_{1}(N) \rightarrow G$ the twisted Alexander
polynomial $\Delta^{\a}_{N,\phi}$ is non--zero. If the subgroup
carried by a connected minimal genus representative of the class
Poincar\'e dual to $\phi$ is separable, then $(N,\phi)$ fibers over
$S^{1}$.
\end{theorem*}

Note that the theorem says in particular that Conjecture 1 holds for all irreducible $N$ with $\pi_1(N)$
surface subgroup separable.
Manifolds with surface subgroup separability  include Seifert
manifolds and, perhaps more importantly, it is conjectured that all
hyperbolic $3$--manifolds satisfy surface subgroup separability (cf. \cite[p.~380]{Th82}).

We point out the somewhat surprising fact that Theorem
\ref{whatever} states that, in the cases under consideration,
Conjecture \ref{conjcha} holds under the apparently much weaker
assumption that the twisted Alexander polynomial is non--zero for
any epimorphism onto a finite group. Furthermore, combined with the
results of \cite{FV06}, this result amounts to the assertion that under the hypotheses of Theorem \ref{whatever}  the set of Seiberg--Witten invariants of all
finite covers of $S^1 \times N$ decide the existence of a symplectic
structure.

Because of their relevance, we quote some corollaries of this
result.

\begin{corollary*} \label{main}
Let $N$ be an irreducible $3$--manifold with vanishing
Thurston norm. If $S^{1} \times N$ is symplectic, then $(N,\phi)$
fibers over $S^1$ for all $\phi\in H^1(N) \setminus \{0\}$.
\end{corollary*}

Recall that  $0$--surgeries along a non--trivial knot are irreducible (cf.
\cite{Ga87}). We can therefore apply Corollary \ref{main} to the
case where the $3$--manifold is obtained as $0$--surgery $N(K)$ of
$S^{3}$ along a knot $K$ of genus $1$. Combined with \cite{BZ67} and
\cite{Ga83} it implies that if $S^{1} \times N(K)$ is symplectic,
then $K$ is a trefoil or the figure--8 knot. This answers in the
affirmative, for the genus $1$ case, Question 7.11 of Kronheimer in
\cite{Kr99} and in particular gives a new proof (see \cite{FV06} for
the original proof) of the fact that if $K$ is the genus--$1$
pretzel knot $(5,-3,5)$, then $S^{1} \times N(K)$ is not symplectic, a
question raised in \cite{Kr98}. Note  that by \cite{Vi03} this corollary
completely characterizes product symplectic manifolds with trivial
canonical class.

Corollary \ref{main} follows from Theorem \ref{whatever} together
with a result of Long and Niblo \cite{LN91}. In Section \ref{sw} we
will also provide a direct and largely self--contained proof based
on, and phrased in terms of, Seiberg-Witten theory for symplectic
$4$-manifolds.

Theorem \ref{whatever} asserts the completeness of twisted Alexander
polynomials for deciding if a pair $(N,\phi)$ satisfying the
hypothesis of this theorem is fibered. In particular, as the
$0$--surgery of $S^{3}$ along $K$ is fibered if and only if $K$ is,
we deduce the following (cf. also \cite{FV06b}).

\begin{corollary*} Twisted Alexander polynomials decide if a knot of genus $1$ is fibered. \end{corollary*}

Building upon Theorem \ref{whatever} and the geometric decomposition
of Haken manifolds we reduce, in Theorem \ref{thm:onlyhyperbolic},
the proof of Conjecture \ref{conjcha} to an appropriate condition of
surface subgroup separability for the hyperbolic components. Note that the condition we require is slightly stronger than surface subgroup separability alone. A corollary of
Theorem \ref{thm:onlyhyperbolic} is worth mentioning:

\begin{corollary*} \label{char} Let $N$ be an irreducible $3$--manifold and let $\phi \in H^{1}(N)$ be a primitive class such that for \textit{any} epimorphism onto a finite group $\a: \pi_{1}(N) \rightarrow G$ the twisted Alexander polynomial $\Delta^{\a}_{N,\phi}$ is non--zero. Denote by $N'$ the union of its Seifert components. Then $(N',\phi|_{N'})$ fibers over $S^{1}$. In particular, Conjecture \ref{conjcha} holds true under the additional assumption that  $N$ is a graph manifold.
\end{corollary*}

The kind of techniques discussed here are applied in \cite{FV07} to study the more general class of symplectic $4$--manifolds admitting a free circle action, obtaining results similar to the ones presented in this paper.

%===================================================
\section{A Seiberg-Witten Proof of Corollary \ref{main}} \label{sw}

Our first goal is to give a proof of Corollary \ref{main} that is as
much as possible self--contained, and that is based on quite
standard results of Seiberg-Witten theory for symplectic
$4$-manifolds and their $3$--dimensional counterpart. In particular,
this proof avoids to deal directly with the somewhat convoluted
definition of Seiberg-Witten invariants for $3$--manifolds with
$b_{1}(N) =1$.

Remember that, for all $\phi \in H^1(N)$, the \emph{Thurston
(semi)norm} of $\phi$ is defined by minimizing the complexity of the
representatives of the class Poincar\'e dual to $\phi$, namely
\[||\phi||_{T}=\min \{ \chi_-(S)\, | \, S \subset N \mbox{ embedded surface dual to }\phi\}.
\]
Here, given a surface $S$ with connected components $S_1\cup\dots
\cup S_k$, we define $\chi_-(S)=\sum_{i=1}^k \max\{-\chi(S_i),0\}$.
By the linearity on rays, this extends to a (semi)norm on
$H^{1}(N;\R)$ (cf. \cite{Th86}).

We have the following straightforward observation.

\begin{lemma} \label{torus}
Let $N$ be an irreducible $3$--manifold, and let $\varphi \in
H^{1}(N;\R)$ be a non--trivial element with $\|\varphi\|_{T} = 0$.
Then $N$ contains a non--separating essential torus $T$. \end{lemma}

\begin{proof}
It is well--known (see e.g. \cite{Th86}) that if the Thurston norm
vanishes for some non--trivial $\varphi\in H^{1}(N;\R)$ then there
exists also a  non--trivial $\phi\in H^{1}(N)$ with $||\phi||_T=0$.
Let $S$ be a (possibly disconnected) embedded surface dual to $\phi$
with $\chi_-(S)=0$. Since cutting $S$ along compressing disks would
increase $\chi_-$ we can assume that each component of $S$ is
incompressible. The hypothesis of irreducibility excludes the case
of spheres. Since $S$ is non--trivial in homology there exists a
connected component $T$ of $S$ that is non--separating. Clearly $T$
satisfies the conditions of the statement.
\end{proof}

We are ready to prove the following theorem, which obviously implies
Corollary \ref{main}.
\begin{theorem} \label{slightly}
Let $N$ be an irreducible $3$-manifold such that $S^{1} \times N$
admits a symplectic form $\w$ whose cohomology class admits
K\"unneth decomposition $[\w] = [dt] \wedge \varphi + \eta$, where
$\varphi \in H^{1}(N;\R)$. If $\|\varphi\|_{T} = 0$ then any $\phi
\in H^1(N;\R) \setminus \{0\}$ can be represented by a closed,
non--degenerate $1$-form; in particular $N$ is a torus bundle and if $\phi$ is a non--trivial integral
class, it can be represented by a fibration. \end{theorem}
\begin{proof}

By Lemma \ref{torus}, we can assume the existence of a
non--separating essential torus $T$ in $N$. In the case when $T$ is
a fiber of a fibration over $S^{1}$, $N$ is a torus bundle. This
means that $N$ is the mapping torus of a self--diffeomorphism $\psi$
of $T^{2}$ classified, up to isotopy, by an element of $SL(2,\Z)$.
The first cohomology group of $N$ is identified with $\Z \oplus
H^{1}(T^{2})^{\psi}$, where $H^{1}(T^{2})^{\psi}$ is the invariant
part of the fiber cohomology (so that $1 \leq b_{1}(N) \leq 3$), and
the Thurston norm vanishes on all of $H^1(N)$. Also, the entire
$H^1(N) \setminus \{0\}$ is composed of fibered classes, and any
nonzero element of the DeRham cohomology is represented by a (unique
up to isotopy) closed, non--degenerate $1$--form (see \cite{Th86}).

We will show now that the case where $T$ is a fiber is the only
possible one. Let us assume, by contradiction, that $T$ is not a
fiber. As $N$ is irreducible and contains a non--separating
essential torus that is not a fiber, it follows from
\cite[Lemma~1~and~Proposition~7]{Ko87} (cf. also \cite{Lu88}) that
the virtual Betti number of $N$ is infinite, and in particular there
exists a finite cover $p : \hat{N} \to N$ with first Betti number
$b_{1}(\hat{N})
> 3$. (Note that this is
excluded in the fibered case: any cover of a torus bundle is itself
a torus bundle, hence the first Betti number, as observed before, is
at most $3$.) The $4$--dimensional Seiberg-Witten polynomial of
$S^{1} \times \hat{N}$ (that has $b_{+}(S^{1} \times \hat{N}) =
b_{1}(\hat{N}) > 1$) coincides, with suitable identification of the
orientations, with the Seiberg-Witten polynomial $SW_{\hat{N}}$ of
$\hat{N}$, an element of $\Z[H^{2}(\hat{N})]$ (see e.g.
\cite{Kr99}). In particular all basic classes $K_i$ are pull-backs,
and we will identify them with elements of $H ^{2}(\hat{N})$.

As there exists a finite covering map $p : S^{1} \times \hat{N} \to
S^{1} \times N$, the manifold $S^{1} \times \hat{N}$ is naturally
endowed with the symplectic form $\hat{\w} := p^{*} \omega$, which
has K\"unneth component $\hat{\varphi} := p^{*}\varphi \in
H^2(\hat{N};\R)$. In particular, the canonical class is a basic
class of $S^{1} \times \hat{N}$, hence is (identified with) an
element of $H^{2}(\hat{N})$ that we denote by $\hat{K}$. We will
exploit now the two main results of Seiberg-Witten theory for
symplectic $4$--manifolds with $b_{+} > 1$, contained in \cite{Ta94}
and \cite{Ta95}. First, the Seiberg-Witten invariant of $\hat{K}$ is
equal to $1$. Second, Taubes' ``more constraints" on the basic classes
$K_{i}$ imply (as $K_{i} \cdot \hat{\w} = K_{i} \cdot
\hat{\varphi}$, the products respectively in $S^{1} \times \hat{N}$
and $\hat{N}$) that
\[ 0 \leq |K_{i} \cdot \hat{\varphi} | \leq \hat{K} \cdot \hat{\varphi},\] and if the latter vanishes, $K_{i} = \hat{K} = 0$.

The $3$--dimensional adjunction inequality for $\hat{N}$ (or, if
preferred, McMullen's inequality relating the Alexander and the
Thurston norm), asserts now that \[ |\hat{K} \cdot \hat{\varphi}|
\leq \|\hat{\varphi} \|_{T} = |\mbox{deg~} p|\|\varphi \|_{T} = 0,\]
where the penultimate equality follows from \cite{Ga87}. This,
together with Taubes' ``more constraints", implies that $\hat{K}$ is
the only basic class and is trivial, which implies in turn that
$SW(\hat{N}) = 1 \in \Z[H^{2}(\hat{N})]$. But it is well-known (see
\cite{Tu01}) that, as $b_{1}(\hat{N}) > 3$, the sum of the
coefficients of $SW(\hat{N})$ (that equals, by \cite{MeTa}, the sum
of coefficients of the Alexander polynomial of $N$) must vanish. We
get therefore a contradiction, which completes the proof.
\end{proof}

Note that Theorem \ref{slightly} covers the case of symplectic $4$--manifolds of the form $S^1 \times N$
having a trivial canonical class. In fact, for these manifolds, the Thurston norm of $N$ must vanish, as discussed in \cite{Vi03}.

%==========================================================================
\section{Twisted Alexander polynomials} \label{sectionufd}

In this section we are going to define (twisted) Alexander polynomial associated to an epimorphism of the fundamental group of a closed $3$--manifold onto a finite group, first
introduced for the case of knots in \cite{L01} (for a broader definition see e.g. \cite{FK06}).

Let $N$ be a closed $3$--manifold and let $\phi: H_{1}(N) \to \Z =
\langle t \rangle$ be a non--trivial homomorphism. We will think of
$\phi$, when useful, as an element of either $Hom(H_1(N),\Z)$ or
$H^{1}(N)$. Through the homomorphism $\phi$, $\pi_{1}(N)$ acts on
$\Z$ by translations. Furthermore let $\a:\pi_{1}(N) \to G$ be an
epimorphism onto a finite group $G$. The composition of $(\a,\phi)$
with the diagonal on $\pi_{1}(N)$ gives an action of $\pi_{1}(N)$ on
$G \times \Z$, which extends to a ring homomorphism from
$\Z[\pi_{1}(N)]$ to the  $\Z[t^{\pm 1}]$--linear endomorphisms of
$\Z[G \times \Z] = \Z[G]\tpm$. This  induces a  left
$\Z[\pi_{1}(N)]$--structure
 on $\Z[G]\tpm$.

Now let $\ti{N}$ be the universal cover of $N$. Note that
$\pi_{1}(N)$ acts on the left on $\ti{N}$ as group of deck
transformation. The chain groups $C_*(\ti{N})$ are in a natural way
right $\Z[\pi_1(N)]$--modules, with the right action on
$C_{*}(\ti{N})$ defined via $\sigma \cdot g := g^{-1}\sigma$, for
$\sigma \in C_{*}(\ti{N})$. We can form by tensor product the chain
complex $C_*(\ti{N})\otimes_{\Z[\pi_1(N)]}\Z[G]\tpm$. Now define
$H_{i}(N;\Z[G]\tpm):=
H_i(C_*(\ti{N})\otimes_{\Z[\pi_1(N)]}\Z[G]\tpm)$, which inherit the
structure of $\Z\tpm$--modules. These module take the name of
twisted Alexander modules.

Our goal is to define an invariant out of $H_{1}(N;\Z[G]\tpm)$.
First note that endowing $N$ with a finite cell structure we can
view
 $C_*(\ti{N})\otimes_{\Z[\pi_1(N)]}\Z[G]\tpm$
as finitely generated $\Z\tpm$--modules. The $\Z\tpm$--module
$H_1(N;\Z[G]\tpm)$ is now a finitely presented and finitely related
$\Z\tpm$--module since $\Z\tpm$ is Noetherian. Therefore
$H_1(N;\Z[G]\tpm)$ has a free $\Z\tpm$--resolution
\[ \Z\tpm^r \xrightarrow{S} \Z\tpm^s \to H_1(N;\Z[G]\tpm) \to 0 \]
of finite $\Z\tpm$--modules. Without loss of generality we can
assume that $r\geq s$.

\begin{defn} \label{def:alex} The \emph{twisted Alexander polynomial} of $(N,\a,\phi)$ is defined
to be the order of the $\Z\tpm$--module $H_1(N;\Z[G]\tpm)$, i.e. the
greatest common divisor of the $s\times s$ minors of the $s\times
r$--matrix $S$. It is denoted by $\Delta_{N,\phi}^{\a}\in \Z\tpm$,
and it is well--defined up to units of $\Z\tpm$.
\end{defn}

Note that this definition only makes sense since $\Z\tpm$ is a UFD.
It is well--known (see e.g. \cite{FV06}) that, up to sign, there is
a unique choice of $\Delta_{N,\phi}^{\a} \in \Z\tpm$ symmetric under
the natural involution of $\Z\tpm$.

If $G$ is the trivial group
we will drop $\a$ from the notation. With these conventions,
$\Delta_{N,\phi} \in \Z[t^{\pm1}]$ is the ordinary $1$--variable
Alexander polynomial associated to $\phi$.

\begin{remark} The $1$-variable twisted Alexander polynomial defined
above can also be described as the  specialization of a
multivariable twisted Alexander polynomial taking values in $\Z[H]$,
where $H$ is the maximal free quotient  of $H_1(N)$. This polynomial, in turn, is
related with the ordinary Alexander polynomial of the $G$--cover
$N_{G}$ of $N$ and then, thanks to \cite{MeTa}, to the
Seiberg-Witten invariants of $S^1 \times N_G$. These observations
constitute the starting point of the connection between Conjecture
\ref{conjcha} and Conjecture \ref{conjfolk}. See \cite{FV06} for
details. \end{remark}

%====================================

%=================================================================

\section{Proof of the Main Theorem} \label{conj}

We will now discuss our main result. Before turning to the
statement, we recall the definition of subgroup separability.

\begin{defn}
Let $\pi$ be a group and $A\subset \pi$ a subgroup. We say that $A$
is \emph{separable} if for any $g \in \pi \setminus A$ there exists
a finite group $G$ and an epimorphism $\a : \pi \to G$ such that
$\a(g) \not\in \a(A)$. A group $\pi$ is called \emph{subgroup
separable} (respectively \emph{surface subgroup separable}) if  any
finitely generated subgroup $A\subset \pi$ (respectively any surface
group $A\subset \pi$) is separable in $\pi$.
\end{defn}

Subgroup separable groups are often also called locally extended
residually finite (LERF).

We are in position now to present our main result.

\begin{theorem} \label{mainthm}
Let $N$ be an irreducible $3$--manifold and let $\phi \in H^{1}(N)$
a primitive class such that for \textit{any} epimorphism onto a
finite group $\a: \pi_{1}(N) \rightarrow G$ the twisted Alexander
polynomial $\Delta^{\a}_{N,\phi}$ is non--zero. If $\phi$ is dual to
a connected incompressible embedded surface $S$ such that $\pi_1(S)$
is  separable in $\pi_1(N)$, then $(N,\phi)$ fibers over $S^{1}$.
\end{theorem}

We point out that the condition that $S$ is connected is not
 restrictive.
 Indeed, McMullen \cite{McM02} showed that if $\phi \in H^1(N)$ is primitive and $\Delta_{N,\phi}\ne 0$, then $\phi$ is dual to a connected incompressible surface.

For the proof of Theorem \ref{mainthm} we will make use of the
following standard result:
\begin{lemma}\label{lem:h0g}
Let $X$ be a connected space, $\a:\pi_1(X)\to G$ a group
homomorphism such that $G/\im(\a)$ is finite. Then
\[ H_0(X;\Z[G])\cong \Z^{|G/Im(\a)|} .\]
\end{lemma}
In fact, the set of components of the (possibly disconnected) finite
cover of $X$ defined by $\a$ gives a basis for  the $\Z$--module
$H_0(X;\Z[G])$ via the Eckmann--Shapiro lemma.

Also, we will need two well--known properties of twisted Alexander
modules.
\begin{lemma}\label{lem:twistedhom}
Let $N$ be a 3--manifold, $\phi\in H^1(N)$ primitive and
$\a:\pi_1(N)\to G$ an epimorphism to a finite group. Then \bn
\item $\Delta_{N,\phi}^\a\ne 0$ if and only if $H_1(N;\zgt)$ is $\zt$--torsion.
\item If $X\subset N$ is a subspace, then
$\rank_{\zt}(H_0(X;\zgt))=0$ if and only if $\phi$ is non--trivial
on $H_1(X)$. Furthermore, if $\phi$ vanishes on $H_1(X)$, then \[
\rank_{\zt}(H_0(X;\zgt)) = \rank_{\Z}(H_0(X;\Z[G]))
=|G|/|\a(\pi_1(X))|.\] \en \end{lemma}

\begin{proof}
The first part is a well--known property of orders. For the second
part note that
 if $\phi$ is non--trivial on $H_1(X)$, then   it follows from Lemma \ref{lem:h0g} applied to
$\a\times \phi:\pi_1(X)\to \Z\times G$ that $H_0(X;\Z[G]\tpm)$ has
finite rank over $\Z$. In particular $H_0(X;\Z[G]\tpm)$ is
$\zt$--torsion. On the other hand, if $\phi$ is trivial on $H_1(X)$,
then $H_0(X;\zgt) =H_0(X;\Z[G])\otimes \zt$, and the lemma follows
from Lemma \ref{lem:h0g}.
\end{proof}

We are now ready to prove  Theorem \ref{mainthm}.

\begin{proof}[Proof of Theorem \ref{mainthm}]
 Let $S\subset N$ be a connected incompressible embedded surface
dual to $\phi$ such that $\pi_1(S)$ is separable in $\pi_1(N)$.  Let
$M : = N \setminus \nu S$, and denote by $\i_{\pm}$ the positive and
negative inclusions of $S$ into $M$. Since $S$ is incompressible,
$(\i_{\pm})_*: \pi_{1}(S) \to \pi_{1}(M)$ is injective, by Dehn's Lemma.
Furthermore it is well--known that this implies that $\pi_{1}(M) \to \pi_{1}(N)$ is
injective as well. We will now show that the hypothesis on $\phi$ imply that either
inclusion induced homomorphism $(\i_{\pm})_*: \pi_{1}(S) \to \pi_{1}(M)$ is in fact an isomorphism.

Pick either inclusion. Denote $A := \pi_{1}(S)$, $B := \pi_{1}(M)$. By the previous observations we can consider
$A$ (via the chosen inclusion) and $B$ as subgroups of $\pi_{1}(N)$. With this notation we have $A\subset B$ and
we have to show that $A=B$.

 Assume, by
contradiction, that there exists an element $g \in B \setminus A$.
Since by assumption $A\subset \pi_1(N)$ is separable, there exist a finite
group $G$ and an epimorphism $\a : \pi_{1}(N) \to G$ such that
$\a(g) \not \in \a(A)$; in particular this implies  $|\a(A)| <
|\a(B)|$.

We will now show  that this contradicts the hypothesis that
$\Delta^{\a}_{N,\phi}$ is non--zero. First note that restricting  the
epimorphism $ \a$ to $A$ and  $ G$ we
define twisted homology modules for $S$ and $M$. These are related with the twisted Alexander modules of $N$ by the following Mayer--Vietoris type exact sequence \[ \ba{cccccccccccccc} \hspace{-0.1cm}&\hspace{-0.1cm}&\hspace{-0.1cm}&\hspace{-0.1cm}\dots \hspace{-0.1cm}&\hspace{-0.1cm}\to\hspace{-0.1cm}&\hspace{-0.1cm} H_1(N;\Z[G]\tpm)\\[0.1cm]
\to \hspace{-0.1cm}&\hspace{-0.1cm}H_0(S;\Z[G])\otimes_\Z \zt \hspace{-0.1cm}&\hspace{-0.1cm}\to
\hspace{-0.1cm}&\hspace{-0.1cm}H_0(M;\Z[G])\otimes_\Z \zt
\hspace{-0.1cm}&\hspace{-0.1cm}\to\hspace{-0.1cm}&\hspace{-0.1cm}
H_0(N;\Z[G]\tpm)\hspace{-0.1cm}&\hspace{-0.1cm}\to\hspace{-0.1cm}&\hspace{-0.1cm} 0. \ea \] (We refer to
\cite[Proposition~3.2]{FK06} for details.) Concerning the terms of the previous sequence, Lemma
\ref{lem:twistedhom} together with the assumption that $\Delta_{N,\phi}^\a\ne 0$ implies that
$H_i(N;\Z[G]\tpm)\otimes_{\zt} \Q(t)=0$ for $i=0,1$.

We are in position now to reach the contradiction. Tensoring the
above exact sequence with $\Q(t)$ we  see  that
\[ \ba{rcl}\rank_{\Z}(H_0(S;\Z[G]))&=&\rank_{\Q(t)}(H_0(S;\Z[G])\otimes_{\Z}\Q(t))\\
&=& \rank_{\Q(t)}(H_0(M;\Z[G])\otimes_\Z \Q(t))=
\rank_{\Z}(H_0(M;\Z[G])).\ea \] It then follows, applying Lemma
\ref{lem:h0g}, that
\[ \frac{|G|}{|\a(A)|}=\rank_\Z(H_0(S;\Z[G]))=\rank_\Z(H_0(M;\Z[G]))=\frac{|G|}{|\a(B)|},\]
which contradicts $|\a(A)|  < |\a(B)|$, hence (reverting to the
standard notation) the maps $\i_{\pm}: \pi_{1}(S) \to \pi_{1}(M)$
are isomorphisms.

Completing the proof is now a standard exercise: note that $\ker\{\pi_1(N)\to \Z\}$ is an infinite amalgamated
product \[ \ker\{\pi_1(N)\to \Z\} =\dots \pi_1(M)*_{\pi_1(S)}*\pi_1(M)*_{\pi_1(S)}*\pi_1(M)\dots,\] where the
inclusion maps are given by $\pi_1(M)\xleftarrow{\i_-}\pi_1(S)$ and $\pi_1(S)\xrightarrow{\i_+}\pi_1(M)$. Since
$\i_{\pm}$ are isomorphisms it follows immediately that $\ker\{\pi_1(N)\to \Z\}\cong \pi_1(M)$, in particular
$\ker\{\pi_1(N)\to \Z\}$ is finitely generated. Since $N$ is irreducible it now follows  from Stallings' theorem
\cite{St62} that $(N,\phi)$ fibers over $S^1$.
\end{proof}

\begin{remark} (1) Scott \cite{Sc78} showed that the fundamental groups of compact surfaces and Seifert manifolds
are subgroup separable. On the other hand it is known (cf.
\cite{BKS87} and \cite{NW01}) that fundamental groups of graph
manifolds and knot complements are in general not subgroup
separable. It is not known whether they are surface subgroup
separable or not. Thurston \cite[p.~380]{Th82} asked whether
fundamental groups of hyperbolic 3--manifolds are subgroup
separable, and various results in this direction are known (see e.g.
\cite{LR05,Gi99}). We refer to \cite{LR05} for more information on
3--manifolds and subgroup separability.
\\ (2) A connected Thurston norm minimizing embedded surface $S$ is incompressible, but the converse is in
general not true. Since there exist separable incompressible Seifert
surfaces for hyperbolic knots which are not of minimal genus (cf.
\cite{AS05}) it might be useful to include, as in the statement
above,  non--Thurston norm minimizing surfaces. Note that in the
case that $(N,\phi)$ fibers over $S^1$ an incompressible surface
dual to $\phi$ is Thurston norm minimizing and unique up to
isotopy.\\
(3) The proof of Theorem \ref{mainthm} carries over to the case that
$N$ has toroidal boundary.
\end{remark}

Whereas subgroup separability is in the general case not completely
understood, the following result of Long and Niblo has particular
relevance for us (see \cite{LN91}).

\begin{theorem}\label{thm:torusseparable} (Long-Niblo)
Let $N$ be a Haken manifold, and $T \subset N$ an embedded
incompressible torus. Then $\pi_{1}(T)$ is separable in
$\pi_{1}(N)$. \end{theorem}

This result has been further generalized by Hamilton, who proved in
\cite{Ha01} that any abelian subgroup is separable in the
fundamental group of Haken manifolds.

The following proposition shows  that
Corollary \ref{main} follows from Theorem \ref{mainthm}:

\begin{proposition}\label{prop:torus}
Let $N$ be an irreducible 3--manifold and $\phi \in H^{1}(N)$ a
primitive class with $||\phi||_T=0$, such that for any epimorphism
onto a finite group $\a : \pi_{1}(N) \to G$ the twisted Alexander
polynomial $\Delta^{\a}_{N,\phi}$ is non--zero. Then $(N,\phi)$
fibers over $S^{1}$.
\end{proposition}

\begin{proof}
 As pointed out after the statement of Theorem \ref{mainthm}, the assumption $\Delta_{N,\phi} \neq 0$ implies that
we can find a connected  Thurston norm minimizing embedded surface $S\subset N$ dual to $\phi$. Clearly $S$ is
an incompressible torus since $N$ is closed, irreducible and $||\phi||_T=0$. The subgroup of $\pi_{1}(N)$
carried by $S$ is separable, hence the statement follows from Theorem \ref{mainthm}.
\end{proof}
%=====================================================================
\section{The JSJ decomposition}

As pointed out in the previous section, $3$--manifolds do not
satisfy, in general, subgroup separability, and it is not clear
whether the weaker condition of surface subgroup separability
required in the hypothesis of Theorem \ref{mainthm} holds or not.
Instead, there is more expectation that some condition of subgroup
separability is satisfied by hyperbolic $3$--manifolds. The goal of
this section is to use the Geometrization Theorem for Haken
manifolds and the results of the previous sections to reduce the
proof of Conjecture \ref{conjcha} to a suitable condition of surface
subgroup separability for hyperbolic manifolds. For manifolds not
already geometric, this is a more direct, and perhaps more realistic
requirement, than the hypothesis of Theorem \ref{mainthm}.

We will start by recalling some standard definitions and results.
(For notation and general results on $3$--manifold topology we refer
to \cite{Bo} and \cite{He76}.)

Let $T_1,\dots,T_s\subset N$ be a family of incompressible embedded
tori. We call $\{T_1,\dots,T_s\}$ a \emph{torus decomposition} if
the (closures of the) components of $N$ cut along $\cup_{j=1}^s T_j$
are either Seifert manifolds or they are simple. (Here simple means
that any incompressible properly embedded torus or annulus is
boundary parallel.)

We call $\{T_1,\dots,T_s\}$ a \emph{JSJ decomposition} if any proper
subfamily fails to satisfy the conditions above. By the work of
Jaco--Shalen and Johannson a JSJ decomposition is unique up to
isotopy. The Geometrization Theorem for Haken manifolds asserts that
the interiors of the simple factors of the decomposition admit a
hyperbolic metric of finite volume.

We are interested in the JSJ decomposition because of the following
theorem.

\begin{theorem}\cite[Theorem~4.2]{EN85} \label{thm:jsjfibered}
Let $N$ be a 3--manifold, $\phi \in H^1(N)$ and $\{T_1,\dots,T_s\}$
a JSJ decomposition. Then $(N,\phi)$ fibers over $S^1$ if and only
if $(N_i,\phi|_{N_i})$ fibers over $S^1$ for every component $N_i$
of $N$ cut along $\cup_{j=1}^s T_i$.
\end{theorem}

This result reduces the problem of fiberability of a $3$--manifold
to the study of its JSJ components. It is natural, within our
approach, to assume that a conjecture similar to Conjecture
\ref{conjcha} holds for manifolds with toroidal boundary. However
note that even if for some nonfibered factor $N_i \subset N$ an
epimorphism $\a : \pi_1(N_i) \to G$ detects nonfiberedness, there is
no reason why that epimorphism should extend to $\pi_1(N)$. This
issue will not cause particular difficulty for the Seifert
components but it is more delicate for the hyperbolic components.
This is analogous to the problem faced in proving residual
finiteness for a Haken manifolds starting from the residual
finiteness of its JSJ components and, in fact, our strategy employs
the pattern of \cite{He87}.

Before we state the next theorem we recall that given a torus $T$
and $\phi \in H^1(T)$, $(T,\phi)$ fibers over $S^1$ if and only if
$\phi \ne 0$.

\begin{theorem}\label{thm:jsjdec}
Let $N$ be an irreducible 3--manifold and $\phi\in H^1(N)$ a
primitive class. Assume that for any epimorphism $\a:\pi_1(N)\to G$
onto a finite group $G$  the twisted Alexander polynomial
$\Delta^{\a}_{N,\phi} \in \Z[t^{\pm 1}]$ is non--zero. Let $T
\subset N$ be an incompressible embedded torus. Then either $\phi|_T
\in H^1(T)$ is non--zero, or $(N,\phi)$ fibers over $S^{1}$ with
fiber $T$.
\end{theorem}

Note that this gives in particular another proof that the examples
in the proof of \cite[Theorem~5.1]{FV06} are not symplectic.

\begin{proof}We start by considering the case where $T$ is non--separating. Assume that $\phi|_T = 0$.
Denote the result of cutting $N$ along $T$ by $M$. As in \cite[Proof~of~Proposition~3.2]{FK06} we get the following
Mayer--Vietoris type exact sequence
\[  \ba{cccccccccccc} &&&&&&H_1(N;\zt)\\[0.1cm]
&\to&  H_0(T;\zt) & \to & H_0(M;\zt)&\to&H_0(N;\zt).\ea \] As $\phi|_T = 0$ and $\Delta_{N,\phi} \ne 0$, Lemma
\ref{lem:twistedhom} implies respectively that $H_0(T;\zt)$ is a non--trivial free $\zt$--module and that
$H_i(N;\zt)$ are $\zt$--torsion modules. Lemma \ref{lem:twistedhom} requires then that $\phi|_M = 0$. This
implies that $T$ is dual to $\phi$ and it follows from Proposition \ref{prop:torus} that $(N,\phi)$ fibers over
$S^1$.

Now assume that $T$ is separating. We will show that $\phi|_T$
cannot be zero. Denote the two components of $N$ cut along $T$ by
$M_1$ and $M_2$. Since $\phi$ is non--zero and the map
$H_1(M_1)\oplus H_1(M_2)\to H_1(N)$ is an epimorphism it follows
that $\phi|_{M_i}$ is non--zero for at least one $i$. Furthermore an
almost identical argument as above shows that if $\phi|_{M_i},
i=1,2$ were both non--zero, then $\phi|_T$ would be non--zero as
well. So we can now assume that $\phi|_{M_i}$ is non--zero for $i=1$
and zero for $i=2$.

Since the kernel of $H_1(T) \to H_1(M_2)$ is nontrivial by Lefschetz
duality, and since $\pi_1(T)=H_1(T)$ it follows that the injective
map $\pi_1(T)\to \pi_1(M_2)$ is not an isomorphism. We can therefore
find  $g\in \pi_1(M_2)\sm \pi_1(T)$. Since $T$ is incompressible we
can view $\pi_1(T)$ and $\pi_1(M_2)$ as subgroups of $\pi_1(N)$. By
Theorem \ref{thm:torusseparable} we can now find an epimorphism
$\a:\pi_1(N)\to G$ onto a finite group $G$ such that $|\a(\pi_1(T))|
< |\a(\pi_1(M_2))|$. In particular
$\rank_\Z(H_0(T;\Z[G]))>\rank_\Z(H_0(M_2;\Z[G]))$. Now consider the
following Mayer--Vietoris type exact sequence
\[  \ba{cccccccccccc} &&&&&&H_1(N;\zgt)\\
&\to&  H_0(T;\zgt) & \to & \bigoplus\limits_{i=1}^{2} H_0(M_i;\zgt)&\to&H_0(N;\zgt).\ea \] It follows
 from Lemma \ref{lem:twistedhom} that, if $\phi|_T = 0$, $H_0(T;\zgt)$ and
$H_0(M_2;\zgt)$ are free $\zt$--modules of ranks
$\rank_\Z(H_0(T;\Z[G]))$ and $\rank_\Z(H_0(M_2;\Z[G]))$.
Furthermore, as $\phi|_{M_1}\ne 0$ and $\Delta_{N,\phi}^\a\ne 0$,
all other modules are $\zt$--torsion modules. However this condition
cannot hold since
$\rank_\Z(H_0(T;\Z[G]))>\rank_\Z(H_0(M_2;\Z[G]))$, hence $\phi|_T \neq 0$.
\end{proof}

In view of Proposition \ref{prop:torus}, we will restrict our
interest to the classes $\phi \in H^1(N)$ with strictly positive
Thurston norm. We can apply Theorem \ref{thm:jsjdec} to the tori of
the JSJ decomposition to prove  the following result.

\begin{proposition}\label{prop:deltacomp}
Let $N$ be an irreducible 3--manifold and  $\phi\in H^1(N)$
primitive with strictly positive Thurston norm. Let
$T_1,\dots,T_s\subset N$ be the JSJ decomposition. Denote the
components of $N$ cut along $\cup_{j=1}^s T_i$ by $N_1,\dots,N_r$,
and let $\phi_i = \phi|_{N_i}$. Assume that for any epimorphism
$\a:\pi_1(N)\to G$ onto a finite group $G$  the twisted Alexander
polynomial $\Delta^{\a}_{N,\phi} \in \Z[t^{\pm 1}]$ is non--zero.
Then for for any epimorphism $\a:\pi_1(N)\to G$ onto a finite group
$G$ and for any $i\in \{1,\dots,r\}$ the twisted Alexander
polynomial $\Delta^{\a}_{N_i,\phi_i} \in \Z[t^{\pm 1}]$ is
non--zero.
\end{proposition}

\begin{proof}
We can apply Theorem \ref{thm:jsjdec} to conclude that $\phi$ is
non--trivial when restricted to $T_i, i=1,\dots,s$. Therefore, for
any epimorphism $\a:\pi_1(N)\to G$ onto a finite group $G$ the
twisted Alexander module $H_1(T_i;\Z[G]\tpm)$ is $\zt$--torsion for
all $i=1,\dots,s$.

Now consider the Mayer--Vietoris type exact sequence
\[ \to \bigoplus_{j=1}^s H_1(T_j;\Z[G]\tpm)\to \bigoplus_{i=1}^r H_1(N_i;\Z[G]\tpm)\to H_1(N;\Z[G]\tpm)\to
\dots.\] Since $H_1(T_j;\Z[G]\tpm), j=1,\dots,s$ and
$H_1(N;\Z[G]\tpm)$ are $\zt$--torsion, it follows that
$H_1(N_i;\Z[G]\tpm), i=1,\dots,r$ are $\zt$--torsion.

This concludes the proof of the proposition.
\end{proof}

\begin{remark} Note that, along the previous lines, it is possible to prove a statement analogous to Proposition \ref{prop:deltacomp} asserting that, if $\Delta^{\a}_{N,\phi}$ is monic, so are the $\Delta^{\a}_{N_i,\phi_i}$. \end{remark}

Theorem \ref{thm:jsjdec} will allow us to control completely the
Seifert components of the JSJ decomposition of a $3$--manifold $N$
that satisfies the hypothesis of the Theorem and, under suitable
assumption of separability, the hyperbolic components.

Before formulating this assumption, we need to recall  some
results and definitions.

First, we will use the classification of incompressible surfaces in
Seifert manifolds. We recall  the following theorem
(\cite[Theorem~VI.34]{Ja80} and \cite[Proposition~1.11]{Hat}).
\begin{theorem}\label{thm:jaco}
Let $N$ be a (compact, orientable) Seifert manifold. If $\S$ is a
connected (orientable)  incompressible surface in $N$, then one of
the following holds: \bn
\item $\S$ is a vertical annulus or torus.
\item $\S$ is a horizontal non--separating surface fibering $N$ as a surface bundle over
$S^1$.
\item $\S$ is a boundary--parallel annulus.
\item $\S$ is a horizontal surface separating $N$ in two twisted $I$--bundles over a compact surface.
\en
\end{theorem}
(Here, a surface in $N$ is called \textit{vertical} (resp.
\textit{horizontal}) if it is the union of fibers (resp. transverse
to all fibers) of some Seifert fibration of $N$.)

Second, observe that given a number $n$ the group $\Z\oplus \Z$ has
precisely one characteristic subgroup of index $n^2$, namely
$n(\Z\oplus \Z)$. Now let $N$ be a 3--manifold with empty or
toroidal boundary. Given a prime $p$ we say that $K\subset \pi_1(N)$
is \emph{$p$--boundary characteristic} if for any component $T$ of
$\partial M$ the group $K\cap \pi_1(T)$ is the characteristic
subgroup of $\pi_1(T)$ of order $p^2$. We denote by $C_p(N)$ the set
of all finite index subgroups of $\pi_1(N)$ which are $p$--boundary
characteristic. (If $N$ has empty boundary, this is just the set of
finite index subgroups.)

We have the following.

\begin{theorem} \label{thm:onlyhyperbolic}
Let $N$ be an irreducible $3$--manifold and let $\phi \in H^{1}(N)$
be a primitive class such that for \textit{any} epimorphism onto a
finite group $\a: \pi_{1}(N) \rightarrow G$ the twisted Alexander
polynomial $\Delta^{\a}_{N,\phi}$ is non--zero. Then the following
hold:
\begin{enumerate} \item  $(N',\phi|_{N'})$ fibers over
$S^{1}$, where $N'$ is the union of the Seifert components.
\item
Assume that  any hyperbolic component $N_i$ satisfies the condition
that, for an incompressible surface $S_i \subset N_i$ Poincar\'e
dual to $\phi|_{N_i}$ and any $g \in \pi_1(N_i) \setminus
\pi_1(S_i)$, there are infinitely many primes $p$ such that there
exists an epimorphism $\pi_1(N_i)\to G$ onto a finite group $G$ with
$\a(g)\not\in \a(\pi_1(S_i))$ and $\ker(\a)\in C_p(N_i)$. Then
$(N,\phi)$ fibers over $S^{1}$. \end{enumerate} \end{theorem}

\begin{proof} If $N$ has a trivial JSJ decomposition, then $N$ is either  a Seifert fibered manifold, or is hyperbolic.
In this  case  the result follows from Theorem \ref{mainthm} because
Seifert manifolds satisfy subgroup separability, respectively
because of the separability assumption of the hypothesis.

We can therefore assume that $N$ has a nontrivial JSJ decomposition
and $\phi$ has strictly positive Thurston norm. In light of Theorem
\ref{thm:jsjfibered}, we want to show that for each component
$N_{i}$, the pair $(N_i,\phi_i)$ is fibered. First note that it
follows from Proposition \ref{prop:deltacomp} that $\phi_i$ is
nontrivial. However, $\phi_i$ is not necessarily primitive, even if
$\phi$ is. Denote by $\varphi_i$ a primitive class in $H^1(N_i)$
with the property that $\phi_i = n \varphi_i$. Clearly
$(N_i,\phi_i)$ fibers over $S^1$ if and only if $(N_i,\varphi_i)$
fibers over $S^1$. Since $\Delta_{N_i,n\varphi_i}(t) =
\Delta_{N_i,\varphi_i}(t^n)$, it follows that
$\Delta_{N_i,\varphi_i} \neq 0$ so that we can find in $N_i$ a
connected minimal genus representative $\Sigma_i$ of the class
Poincar\'e dual to $\varphi_i$.

At this point, we will treat separately Seifert and hyperbolic
components.

Let $N_i$ be a Seifert component. For any component $T$ of $\partial
N_i$, the intersection $\S_{i} \cap T$ is homologically essential,
as $\phi|_{T}$ is a multiple of its Poincar\'e dual, and the former
is nonzero by Theorem \ref{thm:jsjdec}. The knowledge of
incompressible surfaces in Seifert manifolds contained in Theorem
\ref{thm:jaco} will allow us quite easily to show  that $\Sigma_{i}$
is a fiber.

As the intersection of $\S_i$ with the boundary components is
homologically essential, $\S_i$ can satisfy only Case (1) and Case
(2) of Theorem \ref{thm:jaco}, and we claim that also in Case (1)
$\S_i$ fibers $N_i$ over $S^{1}$. This follows by applying
\textit{verbatim} the proof of Theorem \ref{mainthm} to the surface
$\S_i$ in $N_i$, using the condition that
$\Delta^{\a}_{N_i,\varphi_i} \neq 0$ and the fact that the
isomorphic image of the abelian group $\pi_{1}(\S_i) = \Z \subset
\pi_1(N_i) \subset \pi_1(N)$ is \emph{separable in $\pi_1(N)$}, by
the aforementioned result of Hamilton (\cite{Ha01}). Together with
Theorem \ref{thm:jsjfibered} this concludes the proof of (1).

Now let $N_i$ be a hyperbolic component. We write
$N^{c}_i=\cup_{j\ne i}N_j$. It follows from \cite[Lemma~4.1]{He87}
that for all but finitely many prime numbers $p$ there exist $K_j\in
C_p(N_j)$ for all $j\ne i$. It then follows from
\cite[Theorem~2.2]{He87} that in fact for all but finitely many
prime numbers $p$ there exists $K'\in C_p(N_i^{c})$.

Assume, by contradiction, that $\S_i$ is not a fiber of $N_i$. By
the above remark and the separability hypothesis, we can find a
prime number $p$ such that there exists $K'\in C_p(N_i^{c})$ and
such that there exists an epimorphism $\pi_1(N_i)\to G$ onto a
finite group $G$ with $|\a(\pi_1(\S_i)| < |\a(\pi_1(N_i))|$ and
$\ker(\a)\in C_p(N_i)$. Applying again \cite[Theorem~2.2]{He87} we
conclude that there exists a finite index subgroup $K\subset
\pi_1(N)$ such that $K\cap \pi_1(N_i)\subset \ker(\a)$. We can and
will assume that $K$ is normal. Now consider the epimorphism
$\b:\pi_1(N)\to H=\pi_1(N)/K$. Its restriction to $\pi_1(N_i)$ fits
into the following   commutative diagram
\[  \xymatrix{\pi_1(N_i) \ar[r] \ar[dr] & \pi_1(N_i)/(K\cap \pi_1(N_i))\subset \pi_1(N)/K = H \ar[d] \\
& \pi_1(N_i)/\ker(\a)=G. } \]
Since $|\a(\pi_1(\S_i))| <
|\a(\pi_1(N_i))|$ it follows that $|\b(\pi_1(\S_i))| <
|\b(\pi_1(N_i))|$. Following again the argument in the proof of
Theorem \ref{mainthm} we deduce that
$\Delta^{\b}_{N_i,\varphi_i}=0$. But this is a contradiction to
Proposition \ref{prop:deltacomp}.

This shows that $(N_i,\phi_i)$ fibers over $S^1$. Together with (1)
and Theorem \ref{thm:jsjfibered} this concludes the proof of (2).
\end{proof}

If $N$ has no hyperbolic components, we have the following:

\begin{corollary} \label{cor:graph}
If $N$ is a graph manifold (i.e. all components in the JSJ
decomposition are Seifert manifolds), then Conjecture \ref{conjfolk}
holds for $N$.
\end{corollary}

This corollary is particularly significant in light of
\cite{NW01}, that asserts that a graph manifold with positive
$b_1(N)$ satisfies subgroup separability if and only if it is
either Seifert fibered, or a torus bundle over $S^1$.


\begin{thebibliography}{10}
\bibitem[AS05]{AS05}
C. Adams, E. Schoenfeld, {\em Totally Geodesic Seifert Surfaces in
Hyperbolic Knot and Link Complements I.},   Geom. Dedicata  116  (2005), 237--247.
\bibitem[Bo02]{Bo} F. Bonahon, {\em Geometric structures on 3-manifolds}, Handbook of geometric topology,  93--164, North-Holland, Amsterdam, (2002)
\bibitem[BZ67]{BZ67} G. Burde, H. Zieschang, {\em Neuwirthsche Knoten und Fl\"achenabbildungen}, Abh. Math. Sem. Univ. Hamburg  31: 239--246 (1967)
\bibitem[BKS87]{BKS87} R. Burns, A. Karrass, D. Solitar, {\em A note on groups with separable finitely
generated subgroups}, Bull. Austral. Math. Soc. 36 , no. 1: 153--160
(1987)
\bibitem[EN85]{EN85}
D.  Eisenbud, W. Neumann, {\em  Three-dimensional link theory and
invariants of plane curve singularities}, Annals of Mathematics
Studies, 110. Princeton University Press, Princeton, NJ, (1985)
\bibitem[FK06]{FK06}
 S. Friedl, T. Kim, \emph{Thurston norm, fibered manifolds and twisted Alexander
polynomials}, Topology, Vol. 45: 929-953 (2006)
\bibitem[FV06a]{FV06}
S. Friedl, S. Vidussi, {\em Twisted Alexander polynomials and symplectic structures}, Preprint
(2006)
\bibitem[FV06b]{FV06b}
S. Friedl, S. Vidussi, {\em Nontrivial Alexander polynomials of knots and links}, Preprint (2006),
to appear in Bull. Lond. Math. Soc.
\bibitem[FV07]{FV07} S. Friedl, S. Vidussi, {\em Symplectic $4$--manifolds with a free circle action}, Preprint (2007).
\bibitem[Ga83]{Ga83} D. Gabai, {\em Foliations and the
topology of 3--manifolds}, J. Differential Geometry 18, no. 3:
445--503 (1983)
\bibitem[Ga87]{Ga87} D. Gabai, {\em Foliations
and the topology of 3--manifolds. III}, J. Differential Geometry 26,
no. 3: 479--536 (1987)
\bibitem[Gi99]{Gi99}
R. Gitik, {\em Doubles of groups and hyperbolic LERF 3-manifolds},
Ann. of Math. (2) 150, no. 3: 775--806 (1999)
\bibitem[Ha01]{Ha01} E. Hamilton, {\em Abelian Subgroup Separability of Haken $3$--manifolds and Closed Hyperbolic $n$--orbifolds}, Proc. London Math. Soc. 83 no. 3: 626--646 (2001)
\bibitem[Hat]{Hat} A. Hatcher, {\em Basic Topology of 3-Manifolds}, notes available at {\em http://www.math.cornell.edu/\~{} hatcher}
\bibitem[He76]{He76} J. Hempel, {\em $3$-Manifolds},
Ann. of Math. Studies, No. 86. Princeton University Press,
Princeton, N. J. (1976)
\bibitem[He87]{He87} J. Hempel, {\em Residual
finiteness for $3$-manifolds}, Combinatorial group theory and
topology (Alta, Utah, 1984), 379--396, Ann. of Math. Stud., 111,
Princeton Univ. Press, Princeton, NJ (1987)
\bibitem[Ko87]{Ko87}
S. Kojima, {\em Finite covers of $3$-manifolds containing essential
surfaces of Euler characteristic $=0$}, Proc. Amer. Math. Soc. 101,
no. 4: 743--747 (1987)
\bibitem[Kr98]{Kr98} P. Kronheimer, {\em Embedded surfaces and gauge theory in three
and four dimensions},  Surveys in differential geometry, Vol. III
(Cambridge, MA, 1996), 243--298, Int. Press, Boston, MA (1998)
\bibitem[Kr99]{Kr99} P. Kronheimer, {\em Minimal genus in $S\sp 1\times
M\sp 3$},  Invent. Math.  135,  no. 1: 45--61 (1999)
\bibitem[Ja80]{Ja80}
W. Jaco, {\em Lectures on three-manifold topology}, CBMS Regional
Conference Series in Mathematics, 43. American Mathematical Society,
Providence, R.I. (1980)
\bibitem[Li01]{L01} X. S. Lin, {\em Representations of knot groups and twisted
Alexander polynomials}, Acta Math. Sin. (Engl. Ser.)  17,  no. 3:
361--380 (2001)
\bibitem[LN91]{LN91}
D. Long, G. Niblo, {\em Subgroup separability and $3$-manifold
groups}, Math. Z. 207, no. 2: 209--215 (1991)
\bibitem[LR05]{LR05}
D. Long, A. W. Reid, {\em Surface subgroups and subgroup
separability in 3-manifold topology},
 Publicacoes Matematicas do IMPA.  25 $\sp {\rm o}$  Coloquio Brasileiro de Matematica.
 (2005)
\bibitem[Lu88]{Lu88} J. Luecke, {\em Finite covers of $3$--manifolds containing essential tori}, Trans. Amer. Math. Soc. 310: 381--391 (1988)
\bibitem[McC01]{McC01} J. McCarthy, {\em On the asphericity of a symplectic $M^{3} \times S^{1}$}, Proc. Amer. Math. Soc. 129: 257--264 (2001)
\bibitem[McM02]{McM02} C. T. McMullen, {\em The Alexander polynomial of a 3--manifold and the Thurston
norm on cohomology}, Ann. Sci. Ecole Norm. Sup. (4) 35, no. 2:
153--171 (2002)
\bibitem[MeT96]{MeTa} G. Meng, C. H. Taubes, {\em SW = Milnor torsion}, Math. Res. Lett. 3: 661--674 (1996)
\bibitem[NW01]{NW01} G. A. Niblo, D. T. Wise, {\em Subgroup separability, knot groups and graph manifolds}, Proc.
Amer. Math. Soc. 129, no. 3: 685--693 (2001)
\bibitem[Sc78]{Sc78} P. Scott, {\em Subgroups of surface groups are almost
geometric}, J. London Math. Soc. (2) 17, no. 3: 555--565 (1978)
\bibitem[St62]{St62}
J. Stallings, {\em On fibering certain 3--manifolds}, 1962 Topology
of 3--manifolds and related topics (Proc. The Univ. of Georgia
Institute, 1961) pp. 95--100 Prentice-Hall, Englewood Cliffs, N.J.
(1962)
\bibitem[Ta94]{Ta94} C. H. Taubes, {\em The Seiberg-Witten invariants and symplectic forms}, Math. Res. Lett. 1: 809--822 (1994)
\bibitem[Ta95]{Ta95} C. H. Taubes, {\em More constraints on symplectic forms from Seiberg-Witten invariants}, Math. Res. Lett. 2: 9--13 (1995)
\bibitem[Th76]{Th76} W. P. Thurston,
{\em Some simple examples of symplectic manifolds}, Proc. Amer.
Math. Soc. 55 (1976), no. 2, 467--468.
\bibitem[Th82]{Th82} W. P. Thurston,
{\em Three dimensional manifolds, Kleinian groups and hyperbolic
geometry}, Bull. Amer. Math. Soc. 6 (1982)
\bibitem[Th86]{Th86} W. P. Thurston, {\em A norm for the homology of 3--manifolds}, Mem.
Amer. Math. Soc. 339: 99--130 (1986)
\bibitem[Tu01]{Tu01} V. Turaev, {\em Introduction to combinatorial torsions}, Birkh\"auser, Basel, (2001)
\bibitem[Vi99]{Vi99} S. Vidussi, {\em The Alexander norm is smaller than the Thurston norm; a
Seiberg--Witten proof}, Prepublication Ecole Polytechnique 6 (1999)
\bibitem[Vi03]{Vi03} S. Vidussi,
{\em Norms on the cohomology of a 3-manifold and SW theory},
Pacific J. Math.  208,  no. 1: 169--186 (2003)
\end{thebibliography}
\end{document}